\documentclass{amsart}
\usepackage{amsfonts}
\usepackage{amsmath,amscd}
\usepackage{amssymb}
\usepackage{amsthm}
\usepackage{newlfont}

 \newtheorem{thm}{Theorem}[section]
 \newtheorem{cor}[thm]{Corollary}
 \newtheorem{lem}[thm]{Lemma}
 \newtheorem{prop}[thm]{Proposition}
 \theoremstyle{definition}
 
 \theoremstyle{remark}

 \numberwithin{equation}{section}

\begin{document}

\title[]
 {Capability of Nilpotent Lie algebras with small derived Subalgebra}

\author[P. Niroomand]{Peyman Niroomand}
\author{Mohsen Parvizi}
\address{School of Mathematics and Computer Science\\
Damghan University of Basic Sciences, Damghan, Iran}
\email{niroomand@dubs.ac.ir}
\address{Department of Pure Mathematics, Faculty of Mathematical Sciences\\
Ferdowsi University of Mashhad, Mashhad, Iran}
\email{parvizi@math.um.ac.ir}




\date{\today}


\begin{abstract}
In this paper, we classify all capable nilpotent Lie algebras with derived subalgebra of dimension at most 1.
\end{abstract}

\maketitle

\section{Introduction}
Its about seventy years past when P. Hall uttered his famous
problem. "For which group $G$ there exists a group $H$ with $G\cong
H/Z(H)$?" He also noticed that finding such groups are important in
classifying $p$-groups. Following Hall and Senior \cite{ph} a group
$G$ with the above property is called capable. One of the famous
results on this concept is due to Baer \cite{ba} where he classified
all capable groups among the direct sums of cyclic groups and hence
determined all capable finitely generated abelian groups. Capable
groups in the class of extra special $p$-groups were characterized
in \cite{be}. They studied the notion of capability and introduced a
central subgroup denoted by $Z^*(G)$ to be
\[\cap\{\phi(Z(E))\ | \ (E,\phi)~\text{is a  central extension}\},\]
and showed a group $G$ is capable if and only if $Z^*(G)=1$.

Another notion having relation to capability is the exterior square
of groups which was introduced in \cite{br}. Using this concept G.
Ellis \cite{el}, introduced the subgroup $Z^{^{\wedge}}(G)$ to be
the set of all elements $g$ of $G$ for which $g\wedge h=1$ for all
$h\in G$. He could prove $Z^{^{\wedge}}(G)=Z^*(G)$ which is an
interesting result.

Recently several properties of finite $p$-groups has found analogues
results for nilpotent Lie algebras. For instance one can see
\cite{sa} which introduced the notion $Z^*(L)$ for a Lie algebra $L$
similar to what Beyl et. al. \cite{be} introduced for groups. They
also could prove that a necessary and sufficient condition for a Lie
algebra $L$ to be capable is $Z^*(L)=0$.

In this paper we intend to prove some results for Lie algebras.
First of all, we show that the two notions $Z^{^{\wedge}}(L)$ and
$Z^*(L)$ for a nilpotent finite dimensional Lie algebra $L$ are the
same then we give a necessary and sufficient condition for an
abelian finite dimensional Lie algebra to be capable. We also
classify capable Heisenberg Lie algebras and finally using a result
obtained by the first author in a joint paper \cite{ni} we give a
necessary and sufficient conditions for nilpotent finite dimensional
Lie algebras with $\mathrm{dim} (L^2)=1$, to be capable.

\section{preliminaries}

In this section we state some lemmas to use in main results.
Throughout this paper all Lie algebras are finite dimensional,
$A(n)$ and $H(m)$ denote the abelian Lie algebra of dimension $n$
and the Heisenberg Lie algebra of dimension $2m+1$, respectively.

\begin{lem}$($See \cite[Theorem 4.4]{sa}$)$.\label{1}
Let $L$ be a Lie algebra and $N$ be an ideal of $L$. Then
$N\subseteq Z^*(L)$ if and only if the natural map
$M(L)\longrightarrow M(L/N)$ is monomorphism.
\end{lem}

\begin{lem}$($See \cite[Proposition 13]{el2} and \cite[Proposition 4.1 (iii)]{sa}$)$.\label{2}
Let $L$ be a Lie algebra and $N$ be a central ideal of $L$. Then the
following sequences are exact.
\begin{itemize}
\item[(i)] $L\wedge N\longrightarrow L\wedge L\longrightarrow L/N\wedge L/N\longrightarrow
0.$
\item[(ii)] $M(L)\longrightarrow M(L/N)\longrightarrow N\cap L^2\longrightarrow
0.$
\end{itemize}
\end{lem}
\begin{cor}\label{c} $N\subseteq Z^{^\wedge}(L)$ if and only if the natural map $L \wedge L \longrightarrow L/N\wedge
L/N$ is monomorphism.
\end{cor}

\begin{lem}$($See \cite{el3}$)$.\label{3} Let $L$ be a Lie algebra then the following sequence is a central
extension.
\[0\longrightarrow M(L)\longrightarrow L\wedge L\longrightarrow L^2\longrightarrow 0.\]
\end{lem}

\begin{cor}\label{4}
Let $A$ be a finite dimensional abelian Lie algebra. Then
$M(A)\cong A\wedge A.$
\end{cor}

The following lemma describes the Schur multipliers of abelian and
Heisenberg algebras.

\begin{lem}$\mathrm{(See}$ \cite[Lemma 3]{es}, \cite[Example 3]{es2} $\mathrm{and}$ \cite[Theorem 24]{mo}$\mathrm{).}$\label{h}
\begin{itemize}
\item[(i)]$\mathrm{dim}~M(A(n)))=\frac{1}{2}n(n-1)$
\item[(ii)]$\mathrm{dim}~M(H(1))=2$.
\item[(iii)]$\mathrm{dim}~M(H(m))=2m^2-m-1$ for all $m\geq 2$.
\end{itemize}
\end{lem}

\begin{thm}\label{ds}$\mathrm{(See}$ \cite[Theorem 1]{es}$\mathrm{).}$
Let $L_1$ and $L_2$ be finite dimensional Lie algebras. Then
\[\mathrm{dim}~(M(L_1\oplus L_2)) = \mathrm{dim}~(M(L_1)) + \mathrm{dim}~(M(L_2))
+ \mathrm{dim}~(L_1/L_1^2\otimes L_2/L_2^2).\]
\end{thm}

\begin{prop}\label{5} Let $L_1$ and $L_2$ be two Lie algebras. Then
\begin{itemize}
\item[(i)] $(L_1\oplus L_2)\wedge (L_1\oplus L_2)\cong(L_1\wedge L_1)\oplus(L_2
\wedge L_2)\oplus(L_1/L_1^2\otimes L_2/L_2^2).$
\item[(ii)] $Z^{^{\wedge}}(L_1\oplus L_2)\subseteq Z^{^{\wedge}}(L_1)\oplus
Z^{^{\wedge}}(L_2)$
\end{itemize}
\end{prop}
\begin{proof} $(i).$ It is a consequence of \cite[Proposition 8]{el2}.

$(ii)$ It is obtained directly by using $(i)$. \end{proof}

\section{Main Results}
In this section, we classify all nilpotent Lie algebras with
derived subalgebra of dimension at most 1. These Lie algebras are
abelian Lie algebras and nilpotent non-abelian Lie algebras with
derived subalgebra of dimension 1. It was proved in \cite{sa} that a Lie
algebra $L$ is capable if and only if $Z^*(L)=0$. The equality
$Z^*(G)=Z^{^{\wedge}}(G)$ holds for any group $G$. Here, we prove
similar result for Lie algebras and deduce that a Lie algebra $L$ is
capable if and only if $Z^{^{\wedge}}(L)=0$. Also, we state some
lemmas for Heisenberg Lie algebras to use them in main results.

\begin{lem}\label{6}
For any Lie algebra $L$, $Z^*(L)=Z^{^{\wedge}}(L)$
\end{lem}

\begin{proof}
Considering Lemmas \ref{1} and \ref{2} (i), we can see that
$M(L)\longrightarrow M(L/Z^{^{\wedge}}(L))$ is a
monomorphism, so $Z^{^{\wedge}}(L)\subseteq Z^*(L)$. On the other
hand, by Lemmas \ref{1} and \ref{2} $(ii)$, we have
$\mathrm{dim}~M(L/Z^*(L))=
\mathrm{dim}~M(L)+\mathrm{dim}~(L^2\cap Z^*(L))$. But
Lemma \ref{3} shows that \[\mathrm{dim}~(L\wedge
L)=\mathrm{dim}~M(L)+\mathrm{dim}~L^2 ~\text{and}\]\[
\mathrm{dim}~\big(L/Z^*(L)\wedge L/Z^*(L)\big)=\mathrm{dim}
~M(L/Z^*(L))+\mathrm{dim}~(L/Z^*(L))^2.\] Using the
isomorphism $(L/Z^*(L))^2\cong L^2/(Z^*(L)\cap L^2)$, we have
\[\mathrm{dim}~( L\wedge L)=\mathrm{dim}~\big(L/Z^*(L)\wedge
L/Z^*(L)\big),\] hence $Z^*(L)\subseteq Z^{^{\wedge}}(L)$ due to
Corollary \ref{c}.
\end{proof}

\begin{lem}\label{7}
Let $H(m)$ be the Heisenberg Lie algebra. Then
\begin{itemize}
\item[(i)]$H(1)\wedge H(1)\cong A(3)$.
\item[(ii)] $H(m)\wedge H(m)\cong A(2m^2-m)$ for all $m\geq 2$.
\end{itemize}
\end{lem}
\begin{proof}
Since $\mathrm{dim}~ H(m)^2=1$, Lemma \ref{3} follows that $H(m)\wedge
H(m)$ is abelian. Invoking Lemmas \ref{3} and \ref{h}, we should
have $\mathrm{dim}~ (H(1)\wedge H(1))=3$ and $\mathrm{dim}~(
H(m)\wedge H(m))=2m^2-m$ for all $m\geq 2$
\end{proof}

Among the Lie algebras the simplest ones are abelian Lie algebras.
Here we classify all abelian Lie algebras of finite dimension which
are capable.

\begin{thm}\label{t1}
$A(n)$ is capable if and only if $n\geq2$.
\end{thm}

\begin{proof}
Since $M(A(1))=0$, Lemma \ref{1} implies $A(1)$ is not
capable. Now, let $n\geq2$ and $I$ be a $k$-dimensional ideal of
$A(n)$. Then, we have \[\mathrm{dim}~
M(A(n)/I)=\frac{1}{2}(n-k)(n-k-1)~\text{and}~
\mathrm{dim}~ M(A(n))=\frac{1}{2}n(n-1),\] so Lemma
\ref{1} implies $I\subseteq Z^{^{\wedge}}(A(n))$ if and only if
$k=0$. Hence, we should have $Z^{^{\wedge}}(A(n))=0$ and the result
holds.
\end{proof}

Ignoring abelian Lie algebras, Heisenberg Lie algebras are probably
the simplest Lie algebras to work with. The following theorem
classifies all capable Heisenberg Lie algebras.

\begin{thm}\label{t3}
$H(m)$ is capable if and only if $m=1$.
\end{thm}

\begin{proof}
First suppose that $m=1$. Owning to Lemma \ref{7}, we have
$\mathrm{dim}~ (H(1)\wedge H(1))=3$. On the other hand, for any
nonzero ideal of $H(1)$ such as $I$, Corollary \ref{4} and Lemma
\ref{h}(i) follows that
\[\mathrm{dim}~\big((H(1)/I)\wedge (H(1)/I)\big)=1.\] Hence $Z^{^{\wedge}}(H(1))$ contains
no nonzero ideal and must be trivial.

Now, assume that $m\geq 2$ similar to the case $m=1$, by using Lemma
\ref{7}, $\mathrm{dim}~(H(m)\wedge H(m)=2m^2-m$ and
 \[\begin{array}{lcl}\mathrm{dim}~\big(H(m)/{H(m)^2}\wedge H(m)/H(m)^2\big)&=&\mathrm{dim}~
M(H(m)/H(m)^2)\vspace{.3cm}\\&=&2m(2m-1)/2=2m^2-m\end{array}\] which
follows $H(m)^2\subseteq Z^{^{\wedge}}(H(m))$
\end{proof}

The direct sum of an abelian Lie algebra and a Heisenberg Lie
algebra has the derived subalgebra of dimension 1 and its
interesting to know which of them are capable. The following theorem
gives a necessary and sufficient condition for capability of such
Lie algebras.

\begin{thm}
Let $L\cong  H(m)\oplus A(k)$ then $L$ is capable if and only if
$m=1$.
\end{thm}

\begin{proof}
We consider three cases as follows
\begin{itemize}
\item[(i)] $m=k=1$;
\item[(ii)] $m=1$ and $k\geq 2$;
\item[(iii)]  $m\geq 2$.
\end{itemize}
In case $(i)$, $L\cong H(1)\oplus A(1)$ and Proposition \ref{5} and
Theorems \ref{t1}, \ref{t3} follow that
\[Z^{^{\wedge}}(L)\subseteq Z^{^{\wedge}}(H(1))\oplus
Z^{^{\wedge}}(A(1))=A(1).\] But
\[\begin{array}{lcl}\mathrm{dim}~ M(L) &=&\mathrm{dim}~ M(H(1))+\mathrm{dim}~
M(A(1))+\mathrm{dim}~ \big(H(1)/H(1)^2\otimes
A(1)\big)\vspace{.3cm}\\&=&2+0+2=4.\end{array}\] On the other hand,
$\mathrm{dim}~ M(L/A(1))=\mathrm{dim}~ M(H(1))=2$, and so
$A(1)\nsubseteq Z^{^{\wedge}}(L)$ which implies that
$Z^{^{\wedge}}(L)=0$.

In case $(ii)$, Proposition \ref{5} and Theorems \ref{t1}, \ref{t3}
deduce that
\[Z^{^{\wedge}}(L)\subseteq Z^{^{\wedge}}(H(1))\oplus
Z^{^{\wedge}}(A(k))=0,\] as required.

Finally in case $(iii)$, we claim that \[L\wedge L\cong
L/H(m)^2\wedge L/H(m)^2,\] and hence $L$ is not capable.

Since $L/H(m)^2\cong H(m)/H(m)^2\oplus A(k)$, we have
\[\begin{array}{lcl}L/H(m)^2\wedge L/H(m)^2&\cong& (H(m)/H(m)^2\wedge
H(m)/H(m)^2)\vspace{.3cm}\\&\oplus& (A(k)\wedge A(k))\oplus
(H(m)/H(m)^2\otimes A(k)).\end{array}\] Thus \[\mathrm{dim}~
\big(L/H(m)^2\wedge L/H(m)^2\big)=2m(2m-1)/2+n(n-1)/2+2mn.\] On the
other hand, \[\mathrm{dim}~ L\wedge L=\mathrm{dim}~
M(L)+\mathrm{dim}~ H(m)^2+n(n-1)/2+2mn.\] Hence $\mathrm{dim}~
\big(L/H(m)^2\wedge L/H(m)^2\big)=\mathrm{dim}~ L\wedge L$, and the
result holds.
\end{proof}

Now the following theorem by the first author in his joint paper
\cite{ni} shows nothing remains to prove.

\begin{thm}\label{t2}
Let $L$ be an $n$-dimensional nilpotent Lie algebra and
$\mathrm{dim}~ L^2=1$ then $L\cong H(m)\oplus A(n-2m-1)$.
\end{thm}


\begin{thebibliography}{20}
\bibitem{ba} R. Baer, Groups with preassigned central and central quotient group, Trans. Amer. Math.
Soc. 44 (1938) 387-412.
\bibitem{es2} P. Batten and E. Stitzinger,  On covers of Lie algebras, Comm.
Algebra (24)14 (1996), 4301--4317.
\bibitem{es}  P. Batten, K. Moneyhun,  E. Stitzinger.
 On characterizing nilpotent Lie algebras by their multipliers. Comm. Algebra 24(1996) 4319--4330.
\bibitem{be}F. R. Beyl, U. Felgner, and P. Schmid, On groups occurring as center factor groups,
J. Algebra 61 (1979), 161-177.
\bibitem{br}R. Brown,  J.-L Loday, Van Kampen theorems for diagrams of spaces, With an appendix by M. Zisman,
 Topology  26 (1987) 311-335.
\bibitem{el}G. Ellis, Tensor products and q-cross modules, J. London Math. Soc. (2)51
(1995), 241-258.
\bibitem{el2}G. Ellis, A non-abelian tensor product of Lie algebras, Glasg. Math. J. 39 (1991) 101–120.
\bibitem{el3} G. Ellis, Nonabelian exterior products of Lie algebras and an exact sequence in the
homology of Lie algebras, J. Pure Appl. Algebra 46 (1987), 111-115.
\bibitem{ph}M. Hall, Jr. and J. K. Senior, The Groups of Order $2^n$ $(n\leq 6)$, Macmillan Co., New
York, 1964.
\bibitem{mo}
K. Moneyhun, Isoclinisms in Lie algebras, Algebras Groups Geom.
11(1994),  9--22.
\bibitem{ni} Niroomand, P., Russo, F. A note on the Schur multiplier of a nilpotent Lie
algebra. Comm. Algebra, in press.
\bibitem{sa} A.R. Salemkar, V. Alamian, H. Mohammadzadeh, Some properties of the Schur multiplier
and covers of Lie Algebras, Comm. Algebra 36 (2008) 697-707.
\end{thebibliography}
\end{document}